# Modelling the Hidden Flexibility of Clustered Unit Commitment

Germán Morales-España, and Diego A. Tejada-Arango

*Abstract*—This paper proposes a Clustered Unit Commitment (CUC) formulation to accurately model flexibility requirements such as ramping, reserve, and startup/shutdown constraints. The CUC is commonly applied in long-term planning models to approximate the units' operational flexibility in power systems due to its computational advantages. However, the classic CUC intrinsically and hiddenly overestimates the individual unit's flexibility, thus being unable to replicate the result of the individual UC. This paper then present a set of constraints to correctly represent the units' hidden flexibility within the cluster, mainly defined by the individual unit's ramping and startup/shutdown capabilities, including up/down reserves. Different case studies show that the proposed CUC replicates the objective function of the individual UC while solving significantly faster. Therefore, the proposed CUC correctly represents the unit's flexibility within the cluster and could be directly applied to long-term planning models without significantly increasing their computational burden.

*Index Terms*—Unit commitment, clustered UC, ramping constraints, reserves, flexibility, renewable energy sources.

## I. Introduction

MANY power systems worldwide now face a sustained and significant growth of variable renewable energy sources (vRES), such as wind and solar, driven by concerns for the environment, energy security and rising fuel prices. To maintain the supply-demand balance, enough system flexibility must be scheduled in advance to compensate for possible variations in load and renewable production. The day-ahead Unit Commitment (UC) is the short-term planning process that is commonly used to schedule these resources at minimum cost. Recent studies [1], [2] include the UC formulation in long-term planning to correctly represent the operational flexibility of the system in expansion investment decisions. However, this commonly leads to a large-scale mixed-integer programs (MIP), which usually results in prohibitive solving times. To reduce the size and combinatorial complexity of the this problem, units are clustered by technologies [1], [3], and the resulting formulation is commonly called Clustered Unit Commitment (CUC). Despite the computational advantages, CUC overestimates some technical characteristics of the individual units within the cluster. For instance, Meus *et al.* [3] show that, under some specific circumstances, the solution of the individual UC differs from the CUC solution due to an overestimation of startup/shutdown capacities as well as minimum and up/down time limits, and they propose a heuristic method to alleviate this issue. In addition, the overestimation of the ramping flexibility has not been pointed out in the literature, which is becoming more relevant in power systems with a high share of vRES demanding high flexibility. In this paper, we propose a new set of constraints to add into the classic CUC in order to accurately model the units' ramping and reserves flexibility as well as the startup/shutdown capacities without significantly increasing the computational burden.

## II. Clustered UC Formulation

### A. Nomenclature

*Indexes and Sets:*

$g \in \mathcal{G}$  Generating units within a cluster, running from 1 to $G$
$t \in \mathcal{T}$  Hourly periods, running from 1 to $T$ hours

*Individual Unit Parameters:*

$RD$  Ramp-down capability [MW/h]
$RU$  Ramp-up capability [MW/h]
$SD$  Shutdown ramping capability [MW/h]
$SU$  Startup ramping capability [MW/h]
$TD$  Minimum down time [h]
$TU$  Minimum up time [h]

*Continuous Non-negative Variables:*

$\tilde{p}_{gt}$  Power output above minimum output of unit $g$ [MWh]
$\tilde{r}_{gt}^+, \tilde{r}_{gt}^-$  Reserve up/down of unit $g$ [MWh]
$p_t$  Cluster power output above minimum output [MWh]
$r_t^+, r_t^-$  Cluster secondary reserve up/down [MWh]
$\widehat{p}_t$  Total cluster power output [MWh]

*Binary Variables:*

$\tilde{u}_{gt}$  Binary variable which is equal to 1 if the unit $g$ is producing above minimum output and 0 otherwise
$u_t$  Integer variable indicating the number of units producing above minimum output
$y_t$  Integer variable indicating how many units start up
$z_t$  Integer variable indicating how many units shut down

### B. Classic Clustered UC Formulation

This section shows the classic CUC formulation, which is simply a scaled 1-unit formulation. Now the variables $u_t, y_t, z_t$ take integer values $\{1, 2, \ldots, G\}$ instead of binary values [3]. For the sake of brevity, here we show the formulation for one cluster, hence we drop the index for different clusters. The objective function minimizes the total system operational cost (i.e., fixed and variable generation cost, startup/shutdown cost, renewable curtailment cost) and is also subject to system-wide constraints, such as demand balance, transmission limits, and total up/down reserve requirements [4], [5]. The commitment, startup/shutdown logic and the minimum up/down times are guaranteed with [4], [5]

$$u_t - u_{t-1} = y_t - z_t \quad \forall t \tag{1}$$

$$\sum_{i=t-TU+1}^{t} y_i \leq u_t \quad \forall t \in [TU, T] \tag{2}$$

$$\sum_{i=t-TD+1}^{t} z_i \leq G - u_t \quad \forall t \in [TD, T]. \tag{3}$$

The following constraint ensures that the clustered unit operates within its power capacity limits for the case $TU \geq 2$:

$$p_t + r_t^+ \leq \left(\overline{P} - \underline{P}\right) u_t - \left(\overline{P} - SU\right) y_t \\ - \left(\overline{P} - SD\right) z_{t+1} \quad \forall t \quad (4)$$

and when $TU = 1$, the following constraints should be used instead:

$$p_t + r_t^+ \leq \left(\overline{P} - \underline{P}\right) u_t - \left(\overline{P} - SD\right) z_{t+1} \\ - \max(SD - SU, 0) y_t \quad \forall t \quad (5)$$

$$p_t + r_t^+ \leq \left(\overline{P} - \underline{P}\right) u_t - \left(\overline{P} - SU\right) y_t \\ - \max(SU - SD, 0) z_{g,t+1} \quad \forall t. \quad (6)$$

The minimum output and the total energy production are obtained as follows:

$$p_t - r_t^- \geq 0 \quad \forall t \quad (7)$$
$$\widehat{p}_t = \underline{P} u_t + p_t \quad \forall t. \quad (8)$$

It is important to highlight that the set of constraints (1)–(8) is the tightest possible representation (convex hull) for the clustered unit operation [5], where now the vertices of $u_t, y_t, z_t$ are directly scaled from $\{0,1\}$ to $\{0,G\}$. Finally, the ramping limits are also written as a scaled version of the 1-unit constraint:

$$p_t - p_{t-1} + r_t^+ \leq RU \cdot u_t \quad \forall t \quad (9)$$
$$-p_t + p_{t-1} + r_t^- \leq RD \cdot u_{t-1} \quad \forall t. \quad (10)$$

To illustrate how (9) and (10) overestimate the units' ramping capabilities, consider the following example: if 10 units are committed ($u_t = 10$) then the ramp limit is ten times the limit of a single unit. But in case that 9 of the 10 units in the cluster are producing at full capacity at time $t$, then only 1 unit is able to ramp up. Therefore, the real ramping up limit of the units within the cluster is equal to the one of a single unit, not ten times as defined in (9). This is the main drawback of the classic CUC formulation related to the ramping constraints. Furthermore, overestimating the ramping capabilities inherently overestimates reserves, since reserves availability directly depends on the available ramp. In the following section, we propose a set of constraints to be added to the classic CUC formulation to overcome this problem.

### C. Proposed Individual Ramping Constraints for CUC

To overcome the overestimation of ramping and reserve flexibility in the classic CUC, we include the following constraints that guarantee that a single unit cannot increase (decrease) its production and/or up (down) reserve if it is already operating at its maximum (minimum) output capacity:

$$\tilde{u}_{1t} \leq 1 \quad \forall t \quad (11)$$
$$\tilde{u}_{g+1,t} \leq \tilde{u}_{gt} \quad \forall g \in [1, G), t \quad (12)$$
$$\tilde{u}_{Gt} \geq 0 \quad \forall t \quad (13)$$
$$\tilde{p}_{gt} + \tilde{r}_{gt}^+ \leq \left(\overline{P} - \underline{P}\right) \tilde{u}_{gt} \quad \forall g, t \quad (14)$$
$$\tilde{p}_{gt} - \tilde{r}_{gt}^- \geq 0 \quad \forall g, t \quad (15)$$

where (11)-(13) order the commitment of the units, where unit 1 is committed first and successively unit $G$ is committed last.

Enforcing this order tightens the proposed MIP formulation and avoids symmetry, by taking advantage that the committing order does not matter since units within the cluster are assumed to be identical [3]. Constraints (14)-(15) limit the production of a single unit.

The total commitment, up reserve, down reserve, and production (above $\underline{P}$) of the cluster are given by

$$u_t = \sum_{g \in \mathcal{G}} \tilde{u}_{gt} \quad \forall t \quad (16)$$
$$r_t^+ = \sum_{g \in \mathcal{G}} \tilde{r}_{gt}^+ \quad \forall t \quad (17)$$
$$r_t^- = \sum_{g \in \mathcal{G}} \tilde{r}_{gt}^- \quad \forall t \quad (18)$$
$$p_t = \sum_{g \in \mathcal{G}} \tilde{p}_{gt} \quad \forall t. \quad (19)$$

To correctly model the startup and shutdown unit's capabilities, (14) should be replaced by (20), (21) or (22): for the case $TU \geq 2$, capacity limits are ensured with

$$\tilde{p}_{gt} + \tilde{r}_{gt}^+ \leq (SU - \underline{P}) \tilde{u}_{gt} + \left(\overline{P} - SU\right) \tilde{u}_{g,t-1} \quad \forall g, t \quad (20)$$
$$\tilde{p}_{gt} + \tilde{r}_{gt}^+ \leq (SD - \underline{P}) \tilde{u}_{gt} + \left(\overline{P} - SD\right) \tilde{u}_{g,t+1} \quad \forall g, t \in [1,T) \quad (21)$$

and when $TU = 1$, the following constraints should be used:

$$\tilde{p}_{gt} + \tilde{r}_{gt}^+ \leq \left(SU - \overline{P} + SD - \underline{P}\right) \tilde{u}_{gt} + \left(\overline{P} - SU\right) \tilde{u}_{g,t-1} \\ + \left(\overline{P} - SD\right) \tilde{u}_{g,t+1} \quad \forall g, t \in [1,T). \quad (22)$$

It is important to highlight that (20)-(22) overcome the startup/shutdown capability overestimation pointed out in [3].

The ramping limits for individual units are guaranteed with

$$\tilde{p}_t - \tilde{p}_{t-1} + \tilde{r}_{gt}^+ \leq RU \cdot \tilde{u}_{gt} \quad \forall g, t \quad (23)$$
$$-\tilde{p}_t + \tilde{p}_{t-1} + \tilde{r}_{gt}^- \leq RD \cdot \tilde{u}_{g,t-1} \quad \forall g, t. \quad (24)$$

Finally, in the proposed CUC formulation, (15) replaces (dominates) (7), and (23) and (24) replace (9) and (10).

The inequalities (11)-(19) are the tightest possible: First, (11)-(13) describe a convex hull since the resulting matrix constraint is totally unimodular (see proposition 3.2 in [6]). Finally, following Lemma 1 in [7], the new variables $\tilde{p}_{gt}, \tilde{r}_{gt}^+, \tilde{r}_{gt}^-, u_t, p_t$ with (14)-(16) and (19) can be added to the convex hull (11)-(13), and the resulting polytope (11)-(19) is also a convex hull. The final proposed clustered formulation is then given by the two convex hulls (1)–(6), (8) and (11)-(19) together with (20)-(24).

## III. NUMERIC EXPERIMENTS

As case studies, two IEEE systems are used: the New England IEEE 39-bus system and the IEEE 118-bus system. The two systems are scaled by 10, i.e., 10 times the demand, the transmission capacity and the number of generators. Therefore, the scaled IEEE 39-bus system has a total of 90 units, resulting in 9 clusters of 10 units each, and the scaled IEEE 118-bus system has a total of 540 units (54 clusters of 10 units). Both systems data are available online at [8]. All models were solved using GUROBI 8.0 with default options on an Intel-i7

TABLE I
CASE STUDIES RESULTS

| System | Reserve | Result | IUC | CCUC | PCUC-S | PCUC-R | PCUC |
|---|---|---|---|---|---|---|---|
| 39-bus | 10% | O.f. [M$] | 1.0070 | 0.9998 | 1.0051 | 1.0051 | 1.0070 |
| | | O.f. Error | - | 0.72% | 0.20% | 0.20% | 0.00% |
| | | Rtime [s] | 4599 | 4 | 6 | 5 | 15 |
| | 5% | O.f. [M$] | 0.9901 | 0.9826 | 0.9880 | 0.9880 | 0.9901 |
| | | O.f. Error | - | 0.75% | 0.21% | 0.21% | 0.00% |
| | | Rtime [s] | 1218 | 1 | 3 | 4 | 14 |
| 118-bus | 5% | O.f. [M$] | 14.4787 | 14.0853 | 14.2463 | 14.1010 | 14.4789 |
| | | O.f. Error | - | 2.72% | 1.61% | 2.61% | -0.001% |
| | | Rtime [s] | 12543 | 170 | 749 | 388 | 810 |
| | 2.5% | O.f. [M$] | 13.9725 | 13.5540 | 13.7247 | 13.5747 | 13.9724 |
| | | O.f. Error | - | 3.00% | 1.77% | 2.85% | 0.001% |
| | | Rtime [s] | 1924 | 131 | 223 | 377 | 406 |

CPU 3.40 GHz with 16GB of RAM. The MIP problems were solved until they reach a relative optimality gap of $10^{-4}\%$ for the IEEE 39-bus system, and 0.01% for the IEEE 118-bus system.

Table I summarizes the main results for both IEEE systems requiring two different values of up and down reserves (expressed as a percentage of the total demand). Table I compares the objective function (o.f.) and runtimes (Rtime) of the following five UC formulations: IUC, the individual UC [4], [5]; CCUC, the classic CUC [3]; PCUC-S, the proposed CUC without startup/shutdown constraints (20)-(22); PCUC-R, the proposed CUC without ramping constraints (23)-(24); and PCUC, the complete proposed CUC, including both ramping and startup/shutdown constraints of individual units.

In both systems, more reserve requirements imply a higher objective function and longer runtimes in all models. In addition, The classic CUC is expected to be always faster than the proposed CUC because it has a lower number of binary variables, constraints, and non-zero elements. In average, the classic CUC has 10 times less number of binary variables and 6 times less equations than the individual UC, while the proposed CUC has 2.5 times less number of binary variables and almost the same number of equations. When the reserve requirements are higher, the CUC-RS solves 311 times faster than IUC for the IEEE 39-bus and 15 times faster for the IEEE 118-bus system.

By taking the individual UC formulation as a benchmark, we can obtain the objective function error for the CCUC and the proposed CUCs. CUC-RS reduces the objective function error to zero in the IEEE 39-bus system and to 0.001% in the IEEE 118-bus system, notice that this error is lower than the tolerance given to the MIP solver. This is a consequence of the individual ramping constraints and startup/shutdown capacities included in the proposed CUC formulation that succeed in capturing the actual ramping, reserve, and startup/shutdown flexibility of the system. The proposed CUC-R and CUC-S show in the IEEE 39-bus system that both set of constraints are equally relevant for the accuracy in the O.f. error. While in the IEEE 118-bus the ramping constraints are more relevant, since the O.f. error increases more without these constraints. However, the O.f. errors in these simplified versions are still lower than those of the classic CUC.

Finally, in order to diminish the symmetry problem, the IUC was also solved for both cases including a random $\pm1\%$ noise into the variable costs. The results show the same runtime trend as in Table I with an increase in the O.f. error up to 0.05% due to the noise introduced.

IV. CONCLUSION

This paper presented a CUC formulation that better represents the ramping and reserves flexibility of the units. The proposed formulation captures the individual ramping limitation and startup/shutdown capabilities, improving the quality of the solution in comparison to the classic CUC formulation and without significantly increasing the computational burden. The proposed CUC model obtains the same objective function and it solves considerably faster than the individual UC. The proposed formulation can be included in more complex models such as the hybrid model in [3] or capacity expansion models. As future research, the improvement of the unit's individual minimum up/down time constraints is one of the challenges to overcome, since the classic CUC fails to accurately represent them. Although our results initially suggest that these constraints have a much lower impact than the ramping limitation and startup/shutdown capabilities addressed in this paper.